\documentclass[a4paper,12pt]{article}
\usepackage{amsfonts}
\usepackage{latexsym}
\usepackage{amsmath}
\usepackage{MnSymbol}
\usepackage{indentfirst}
\usepackage{sectsty}
\newcommand{\ms}{\scriptscriptstyle} 

\begingroup

\endgroup

\begin{document}

\noindent{{\Large{\textbf{Rough sets and three-valued structures}}}
\footnote{This article is published in \textit{Logic at Work}, chapter 33, Or{\l}owska E. (ed.), Essays Dedicated to the Memory of Helena Rasiowa, Physica-Verlag, Heidelberg, 1999, 596--603.}

\

\noindent \textbf{Luisa ITURRIOZ}

\

\noindent \textit{Universit\'e de Lyon, Universit\'e Claude Bernard Lyon 1, \\
Laboratoire de Math\'ematiques Discr\`etes,
F-69622 Villeurbanne cedex, France.} \\
\noindent E-mail: luisa.iturrioz@math.univ-lyon1.fr}

\

\hfill{To the memory of Helena Rasiowa}

\

\noindent \textbf{Abstract:} In recent years, many papers have been published showing relationships between rough sets and some lattice theoretical structures. We present here some strong relations between rough sets and three-valued\\
{\L}ukasiewicz algebras.

\medskip

\noindent \textbf{Key Words:} 

\noindent {\small Information systems, rough sets, monadic Boolean algebras, {\L}ukasiewicz and Post algebras of order 3}

\

\noindent \textbf{\S{1}. Introduction}

\bigskip

The concept of a rough set was introduced by Pawlak (1981, 1982) as a tool to approximate a set by a pair of sets, called the lower and the upper approximation of this set.

In the last few years many papers have been published showing relations between the concepts of an approximation space and rough sets on the one hand, and lattice theoretical structures (distributive lattices, Stone algebras, regular double Stone algebras, semi-simple Nelson algebras, etc$\ldots$) on the other hand.

Also, a lot of fundamental research related to multiple-valued logics has been developed by Helena Rasiowa. This paper is an attempt to state relations between rough sets and finite-valued ({\L}ukasiewicz, Post) algebras. In this way, we hope to contribute to a better understanding of the rough set notion. 

For the sake of clarity we collect in \S{2} the definitions and known results which will be used throughout the paper. In \S{3} we exhibit a strong relation between rough sets and three-valued {\L}ukasiewicz algebras. This relation was inspired by the Moisil construction [1940] of centered three-valued {\L}ukasiewicz algebras via Boolean algebras. In \S{4}, we give an answer to the problem concerning the extensions of membership functions on rough sets by applying the Monteiro construction [1967] of three-valued {\L}ukasiewicz algebras using Monadic Boolean algebras.

\

\noindent \textbf {\S{2}. Preliminaries}

\bigskip

In this section we recall some basic notions related to equivalence relations and monadic Boolean algebras [1], [2], information systems and rough sets [13]. 

Let $Ob$ be a nonempty set (set of objects) and $R$ an equivalence relation on $Ob$. Let$R^{\ast}$ be the family of all equivalence classes of $R$, i.e. $R^{\ast} = \{\vert x \vert : x \in Ob\}$. This family is a partition of $Ob$. 

It is well known that on the Boolean algebra ${\cal B} = ({\cal P}(Ob), \cap, \cup, \neg,\emptyset, Ob)$ where ${\cal P}(Ob)$ denotes the powerset of $Ob$, and $\neg A = Ob - A$, the equivalence relation $R$ induces a unary operator $M$ in the following way, for $A \subseteq Ob$: 
\[
MA = \bigcup \{\vert x \vert \in R^{\ast} : x \in A\};
\]

\noindent which is equivalent to

\[
MA = \bigcup \{ \vert x \vert \in R^{\ast} : \vert x \vert \cap A \neq \emptyset\}.
\]

\noindent By definition we have $M(\emptyset) = \emptyset$ and $A \subseteq MA$. It is well known (see for example [1], [16]) that $M$ also satisfies the condition $M(A \cap MB) = MA \cap MB$, for all $A,B \in {\cal P}(Ob)$. For the sake of clarity we recall the proof of this equality. Let $z \in M(A \cap MB)$ then there exists $x \in A \cap MB$ such that $z \in \vert x \vert$. Since $x \in A \cap MB$ we infer $x \in A$ and there exists $y \in B$ such that $x \in\vert y \vert$. Thus $\vert x \vert = \vert y \vert$. Therefore $z \in MA \cap MB$. Conversely, let $z \in MA \cap MB$ then there exist $x \in A$ and $y \in B$ such that $z \in \vert x \vert$ and $z \in \vert y \vert$, so $\vert x \vert = \vert y \vert$. Therefore $x \in A \cap MB$ and $z \in M(A \cap MB)$.

\

We conclude that $M$ is a \textbf {monadic operator} on the Boolean algebra ${\cal B}$. Recall that in [2] a \textbf {monadic operator} (or quantifier, or S5 operator) $M$ on a Boolean algebra ${\cal B} = (B, \wedge, \vee, \neg, 0,1)$ is a map $M : B \rightarrow B$ satisfying the following conditions:
\[
\begin{array}{lllll}&(M0)&M0= 0\\
&(M1)&a \wedge Ma = a\\
&(M2)&M(a \wedge Mb) = Ma \wedge Mb.\\
\end{array}
\]

\noindent The system ${\cal B} = (B, \wedge, \vee, \neg, 0,1, M)$ of type $(2,2,1,0,0,1)$ is called a {\bf monadic Boolean algebra}. For equivalent definitions see [1].

The notion of monadic Boolean algebra has been introduced by Halmos [2] in order to give a systematic algebraic study of the classical monadic functional calculus; the operation $M$ is called the existential quantifier. As usual the universal quantifier is defined by $Lx = \neg M \neg x$. In ${\cal P}(Ob)$, the operator $L$ is defined by
\[
LA = \bigcup\{\vert x \vert \in R^{\ast} : \vert x \vert \subseteq A\},\ \ for\ A \subseteq Ob.
\]

\noindent Monadic Boolean algebras have many interesting algebraic properties [8]. In [2], [3] it has been shown that the set of closed elements (fixed points) $M(B) = \{x \in B : Mx = x\} =\{x \in B : Lx = x\}$ is a monadic Boolean subalgebra of ${\cal B}$.

\

An \textbf {information system} in the sense of Pawlak [13] is a system 
\[
(Ob, Att, \{Val_{a} : a \in Att\}, f)
\]

\noindent where $Ob$ is a nonempty (finite) set called the universe of objects, $Att$ is a nonempty finite set of attributes, each $Val_{a}$ is a nonempty set of values of attribute $a$, and $f$ is a function $f : Ob\ \times Att \rightarrow Val$, where $Val = \bigcup_{a \in Att} Val_{a}$. In this way, for every $x \in Ob$ and $a \in Att$ we have that $f(x,a) = a(x) \in Val_{a}$.

\

An \textbf {equivalence relation} $R$ on $Ob$, called the {\bf indiscernibility relation}, can be defined in the following way: 

for $x, y \in Ob,\ \ x R y$ if and only if $f(x,a) = f(y,a)$, for every $a \in Att$. The system $(Ob, R)$ is called an \textbf{approximation space.}

\

It follows by the construction above that this equivalence relation generates a monadic operator $M$ and its dual $L$ on the Boolean algebra\\
 ${\cal B} = ({\cal P}(Ob), \cap, \cup, \neg, \emptyset, Ob)$. For each $X \subseteq Ob$ we have elements $LX$ and $MX$ of the monadic Boolean subalgebra $M({\cal P}(Ob))$ with $LX \subseteq X \subseteq MX$.

\

\noindent \textbf {\S{3}. Rough sets and three-valued {\L}ukasiewicz algebras}

\bigskip

Recall that a \textbf{three-valued {\L}ukasiewicz algebra} [5], [9] is an algebra\\
$(L, \wedge, \vee, \sim,\bigtriangledown, 1)$ of type $(2,2,1,1,0)$, where $(A,\wedge, \vee, \sim, 1)$ is a De Morgan algebra and $\bigtriangledown$ is a unary operator (the possibility operator) satisfying the following conditions: 

$\sim x \vee\bigtriangledown x =1 ;\ x \wedge \sim x =\ \sim x \wedge \bigtriangledown x ;\ \bigtriangledown(x \wedge y) = \bigtriangledown x \wedge \bigtriangledown y$.

The three-valued {\L}ukasiewicz logic [4] has an algebraic interpretation in three-valued {\L}ukasiewicz algebras. The operator $\bigtriangledown$ is an additive-multipli\-cative closure operator such that the set $\bigtriangledown(L) =\{x \in L : \bigtriangledown x = x\}$ of invariant elements is the Boolean subalgebra of complemented elements. The necessity operator is defined by $\bigtriangleup x =\ \sim \bigtriangledown \sim x$.

Moisil [5] has proved that {\L}ukasiewicz algebras satisfy the following ``determination principle": If $\bigtriangledown x = \bigtriangledown y$ and $\bigtriangleup x = \bigtriangleup y$ then $x=y$.

A \textbf{centered three-valued {\L}ukasiewicz algebra}, or a \textbf{Post algebra of order 3}, is a three-valued {\L}ukasiewicz algebra with a \textbf{center}, that is an element $c$ of $L$ such that $\sim c = c$. The center of $L$ (if it exists) is unique, and $x = (\bigtriangleup x \vee c) \wedge \bigtriangledown x$, for all $x \in L$ [7]. Other equivalent definitions of Post algebras of order 3 can be found in [14].

A \textbf{rough set} of the approximation space $(Ob, R)$[13] is a pair $(LX, MX)$ where $X \subseteq Ob$. Let $B^{\ast}$ be the collection of all rough sets of $(Ob, R)$. We will define on $B^{\ast}$ some algebraic structure.

Since $LX$ and $MX$ are elements of the monadic Boolean algebra $M({\cal P}(Ob))$ of closed elements of $M$ and $LX \subseteq X \subseteq MX$, we consider the following operations on $B^{\ast}$:
\begin{align*}(LX, MX) \wedge (LY, MY)&= (L(LX \cap LY), M(MX \cap MY))\\
 &= (LX \cap LY, MX \cap MY) = (LU, MU) \quad (\ast) [12]\\
(LX, MX) \vee (LY, MY)&= (L(LX \cup LY), M(MX \cup MY)) \\
&= (LX \cup LY, MX \cup MY) = (LV, MV) \quad (\ast \ast) [12]\\
\sim(LX, MX)&= (L \neg MX, M \neg LX)\\
&= (L \neg X, M \neg X)\\
\bigtriangledown(LX, MX)&=(LMX, MMX) = (MX, MX)\\
0 = (\emptyset, \emptyset)  \ ; \ \ \ \ 1 = (Ob, Ob)\\
\end{align*}

The right side equalities above are in $B^{\ast}$ [12] because the system\\
$(M({\cal P}(Ob)),\cap, \cup, \neg, \emptyset, Ob, M)$ is a monadic Boolean subalgebra of ${\cal B}$.

\

We show the following result.

\medskip
\noindent \textbf{Theorem 3.1}.

For every approximation space $(Ob, R)$, the system $(B^{\ast}, \wedge, \vee, \sim,\bigtriangledown, 1)$ is a three-valued {\L}ukasiewicz algebra.

\medskip

\noindent \textbf{Proof}. The proof is straightforward.

\medskip

\noindent \textbf{Remark}. $(B^{\ast}, \wedge, \vee, \sim, \bigtriangledown, 0, 1)$ is a three-valued {\L}ukasiewicz subalgebra of the \textbf{centered} three-valued {\L}ukasiewicz algebra that we can obtain by applying the Moisil construction [5, p. 450], [7, p. 200] from the Boolean algebra $M({\cal P}(Ob))$. Recall that in the Moisil construction, we consider the set of all pairs $(b_{1}, b_{2})$ of elements of a Boolean algebra with $b_{1} \leq b_{2}$. The center is the element $c = (\emptyset, Ob)$.

\medskip

The construction above shows that all the results known in the theory of three-valued {\L}ukasiewicz algebras can be applied to rough sets. For example, we know that every three-valued {\L}ukasiewicz algebra is a Heyting algebra [6]; the intuitionistic implication being defined by the equality:

$x \Rightarrow y =\ \sim \bigtriangledown x \vee y \vee(\bigtriangledown \sim x \wedge \bigtriangledown y).$

Also, every three-valued {\L}ukasiewicz algebra is a Kleene algebra (i.e. the condition $x\ \wedge \sim x \leq y\ \vee \sim y$ is satisfied) [9], [11], a Stone algebra [9], [11], a regular double Stone algebra [15], a semi-simple Nelson algebra [11], etc...

\medskip

We are interested here in showing a converse result, namely, that every three-valued {\L}ukasiewicz algebra can be represented as an algebra of rough subsets of an approximation space $(Ob, R)$.

First we present some definitions and results which hold in the theory of {\L}ukasiewicz algebra of order 3.

Let $A$ be a three-valued {\L}ukasiewicz algebra and let $Ob$ be the set of all prime filters in $A$, ordered by inclusion. We consider the Bia{\l}ynicki-Birula and Rasiowa order reversing involution $g : Ob \rightarrow Ob$ defined in the following way [14, p. 45]: for any $P \in Ob,\ g(P) = \neg \sim P$ where $\sim P =\{\sim p : p \in P\}$ and $\neg$ is the set theoretical complement.

The set $Ob$, ordered by inclusion, is the disjoint union of chains of one or two elements [11]. Since $A$ satisfies the Kleene law then for each $P \in Ob$, prime filters $P$ and $g(P)$ are comparable [Rasiowa 1958], [11, p. 45].

If $P, Q \in Ob$ then we define $P R_{\ms Ob} Q$ if and only if $P$ and $Q$ are comparable, i.e. if they are in the same chain. $R_{\ms Ob}$ is an equivalence relation on $Ob$ such that if $P R_{\ms Ob} Q$ then $g(P) R_{\ms Ob} g(Q)$.

Let us note the following results.

\medskip

\noindent \textbf{Lemma 3.1.} If $g(P) \subseteq P$ and $\bigtriangledown x \in P$ then $x \in P$.

\

\noindent \textbf{Proof.} If $x \not \in P$ then $\sim x \not \in\ \sim P$ so $\sim x \in g(P)\subseteq P$. Since $P$ is a filter, $\bigtriangledown x \in P$ and $x\ \wedge \sim x =\ \sim x \wedge\bigtriangledown x$ it follows that $x\ \wedge \sim x \in P$, so $x \in P$, a contradiction.

\

\noindent \textbf{Lemma 3.2.} If $g(P) \subseteq P$ and $x \in g(P)$ then $\bigtriangleup x \in g(P)$.

\

\noindent \textbf{Proof.} If $\bigtriangleup x \not \in g(P)$ then $\sim \bigtriangleup x = \bigtriangledown \sim x \in g(P) \subseteq P$ because $g(P)$ is a prime filter and $1 = \bigtriangleup x\ \vee \sim \bigtriangleup x \in g(P)$. Hence $x \wedge \bigtriangledown \sim x = x \wedge \sim x \in P$ and $\sim x \in P$. Or since $x \in g(P) = \neg \sim P$ we get $\sim  x \not \in P$, a contradiction.

\

We will show that collections of rough sets of an approximation space are typical examples of three-valued {\L}ukasiewicz algebras in the sense indicated in the next theorem.

\

\noindent \textbf{Theorem 3.2. Representation Theorem}. Every three-valued {\L}ukasiewicz algebra $A$ can be represented as an algebra of rough subsets of an approximation space $(Ob, R_{\ms Ob})$.

\

\noindent \textbf{Proof}. Let $Ob$ be the set of all prime filters in $A$, ordered by inclusion and $R_{\ms Ob}$ the equivalence relation defined above. We consider the monadic Boolean algebra ${\cal B} = ({\cal P}(Ob), \cap, \cup, \neg, \emptyset, Ob, M)$, where $MX = \bigcup \{\vert P \vert\in R^{\ast}_{\ms Ob} : P \in X\}$, for $X \subseteq Ob$.

Following Stone, for every $x \in A$ we define the map $s : A \rightarrow {\cal P}(Ob)$ as follows: $s(x) = \{P \in Ob : x \in P\}$. The map $s$ is a one-one $(0,1)$-lattice homomorphism.

Let $B^{\ast}$ be the collection of pairs $(Ls(x), Ms(x))$ with operations defined in the above way. The system $(B^{\ast}, \cap, \cup, \sim, \bigtriangledown, \emptyset, Ob)$ is a three-valued\\
{\L}ukasiewicz algebra.

We consider the map $h : A \rightarrow B^{\ast}$ defined as follows: $h(x) = (Ls(x), Ms(x))$.

We are going to show that $h$ is a three-valued {\L}ukasiewicz isomorphism. 

\noindent We need the following results.

\medskip

\noindent $(1)$ $Ms(\bigtriangledown x) = s(\bigtriangledown x)$.

By definition of $M$ we have $s(\bigtriangledown x) \subseteq Ms(\bigtriangledown x)$. On the other hand let $P \in Ms(\bigtriangledown x)$, so there exists $Q \in Ob$ such that $P R_{\ms Ob}Q$ and $Q \in s(\bigtriangledown x)$, i.e. $\bigtriangledown x \in Q$; if $Q \subseteq P$ then $\bigtriangledown x \in P$ ; if $P \subset Q$ and $\bigtriangledown x \not \in P$ then since $\sim \bigtriangledown x \vee \bigtriangledown x = 1 \in P$ and $P$ is prime we obtain $\sim \bigtriangledown x \in P \subset Q$ so $\sim \bigtriangledown x \wedge  \bigtriangledown x = 0 \in Q$, a contradiction. In both cases $P \in s(\bigtriangledown x)$.

\medskip

\noindent $(2)$ $Ms(x) = s(\bigtriangledown x)$.

By $x \leq \bigtriangledown x$ then $s(x) \subseteq s(\bigtriangledown x)$ and since $M$ is order preserving we obtain (i) $Ms(x) \subseteq Ms(\bigtriangledown x) = s(\bigtriangledown x)$ by (1). To show that (ii) $s(\bigtriangledown x) \subseteq Ms(x)$ we suppose $P \in s(\bigtriangledown x)$, i.e. $\bigtriangledown x \in P$. If $g(P) \subseteq P$ then by lemma 1 we obtain $x \in P$, so $P \in s(x) \subseteq Ms(x)$. If $P \subseteq g(P)$ then $gg(P) \subseteq g(P)$ and $\bigtriangledown x \in g(P)$; by lemma 1, $x \in g(P)$ so $g(P) \in s(x)$ and $P \in Ms(x)$.

\medskip

\noindent $(3)$ $Ls(\bigtriangleup x)  = s(\bigtriangleup x)$.

By definition of $L$ we have $Ls(\bigtriangleup x) \subseteq s(\bigtriangleup x)$. On the other hand let $P \in s(\bigtriangleup x)$ so $\bigtriangleup x \in P$. If $g(P) \subseteq P$ then $\bigtriangleup x \in g(P)$ (in fact, if $\bigtriangleup x \not \in g(P)$ then $\sim \bigtriangleup x \in P$ and since $\bigtriangleup x \in P$ we have $0 \in P$, a contradiction). If $P \subset g(P)$ then$\bigtriangleup x \in g(P)$. In both cases $\vert P \vert \subseteq s(\bigtriangleup x)$ and $P \in Ls(\bigtriangleup x)$.

\medskip

\noindent $(4)$ $Ls(x) = s(\bigtriangleup x)$.

By $\bigtriangleup x \leq x$ it follows that $s(\bigtriangleup x) \leq s(x)$ and since $L$ is order preserving we obtain $Ls(\bigtriangleup x) \subseteq Ls(x)$ and by $(3)$, $(i)$ $s(\bigtriangleup x) \subseteq Ls(x)$. To show that $(ii)$ $Ls(x) \subseteq s(\bigtriangleup x)$ we suppose $P \in Ls(x)$ then $\vert P \vert \subseteq s(x)$; this implies that $x \in P$ and $x \in g(P)$. If $g(P) \subseteq P$ then by lemma 2, $\bigtriangleup x \in g(P) \subseteq P$. If $P \subseteq g(P)$ then by lemma 2 again, $\bigtriangleup x \in P$. In both cases $P \in s(\bigtriangleup x)$ as required.

\medskip

\noindent $(5)$ $s(\sim \bigtriangledown x) = \neg s(\bigtriangledown x) = \neg Ms(x)$.

In fact it is a consequence of the following equivalent conditions:

\noindent $P \in s(\sim \bigtriangledown x) \Leftrightarrow\ \sim \bigtriangledown x \in P \Leftrightarrow\bigtriangledown x \not \in  P \Leftrightarrow P \not \in s(\bigtriangledown x) \Leftrightarrow P\in \neg s(\bigtriangledown x)$

\medskip

\noindent $(6)$ $s(\sim \bigtriangleup x) = \neg s(\bigtriangleup x) = \neg Ls(x)$

It is a consequence of the following equivalent conditions:

\noindent $P \in s(\sim \bigtriangleup x) \Leftrightarrow\  \sim \bigtriangleup x \in P \Leftrightarrow\bigtriangleup x \not \in P \Leftrightarrow P \not \in s(\bigtriangleup x) \Leftrightarrow P \in \neg s(\bigtriangleup x)$

\medskip

Using the results $(1)$-$(6)$ it is straightforward to show that $h$ is a three-valued {\L}ukasiewicz homomorphism. The map $h$ is one-one. Indeed, suppose that $(s(\bigtriangleup x), s(\bigtriangledown x)) = (s(\bigtriangleup y), s(\bigtriangledown y))$ then $s(\bigtriangleup x) = s(\bigtriangleup y)$ and $s(\bigtriangledown x) = s(\bigtriangledown y)$. Since $s$ is one-one we have $\bigtriangleup x = \bigtriangleup y$ and $\bigtriangledown x= \bigtriangledown y$. By the Moisil determination principle we conclude $x=y$.

The proof of the representation theorem is now complete.

\

\noindent\textbf {\S{4}. Rough sets and membership functions}

\bigskip

Let ${\cal B} = ({\cal P}(Ob), \cap, \cup, \neg, \emptyset, Ob, M)$ be the monadic Boolean algebra generated by the indiscernibility relation $R$ on an information system. Following A. Monteiro[10] we define $A \equiv B (mod R)$ if and only if $LA = LB$ and $MA = MB$ and two new operations on ${\cal P}(Ob)$ by 
\begin{align*}
A \capdot B &= MA \cap B \cap (A \cup M \neg B)\\
A \uplus B& = LA \cup B \cup (A \cap L \neg B)
\end{align*}

The relation $\equiv$ is a congruence on $\capdot, \uplus, \neg, M$. Let ${\cal P}(Ob) / \equiv$ be the set of all equivalence classes and $\vert A \vert$ the equivalence class containing $A$. If we consider ${\cal P}(Ob)/ \equiv$ algebrized in the natural way, i.e.\ $1 = \vert Ob \vert, \vert A \vert \capdot \vert B \vert = \vert A \capdot B \vert, \vert A \vert \uplus \vert B \vert = \vert A\uplus B \vert, \sim \vert A \vert = \vert \neg A \vert, \bigtriangledown \vert A \vert = \vert MA\vert$ we have that the system $({\cal P}(Ob) / \equiv, \capdot, \uplus, \sim, \bigtriangledown, 1)$ is a three-valued {\L}ukasiewicz algebra [10], [12].

In particular the following equalities are satisfied [12]:
\begin{align*}
M(A \capdot B)&=MA \cap MB, &M(A \uplus B) &= MA \cup MB,&\\
L(A \capdot B)&=LA \cap LB, &L(A \uplus B) &= LA \cup LB.
\end{align*}

Let $A \subseteq Ob$. We define a membership function associated to $A$ in the following way [Pawlak,1985, [13]]:
\[
 \mu_{A}(x) = \left\{ 
 \begin{array}{llll} 
 1 &\text{iff}& x \in LA\\
 1/2 &\text{iff}& x \in MA \cap \neg LA\\
 0 &\text{iff}& x \in \neg MA\\
\end{array} 
\right.
\]

\noindent \textbf{Proposition 4.1.} The membership function $\mu_{A} : Ob \rightarrow \{0,1/2, 1\}$ can be extended to the operations $\capdot, \uplus$ and $\neg$ between sets.

\

\noindent \textbf{Proof.} We note $\lambda$ this extension. By the way of example we show the validity of the extension for $A \uplus B$. It is a consequence of the following equivalent conditions:

$\lambda_{A \uplus B}(x) = 1 \Leftrightarrow x \in L(A \uplus B) = LA \cup LB\Leftrightarrow x \in LA$ or $x \in LB \Leftrightarrow \mu_{A}(x) =1$ or $\mu_{B}(x) = 1 \Leftrightarrow max(\mu_{A}(x), \mu_{B}(x)) = 1$ 

$\lambda_{A \uplus B}(x) = 1/2 \Leftrightarrow x \in M(A \uplus B) \cap \neg L(A \uplus B) = (MA \cup MB) \cap \neg (LA \cup LB) = (MA \cup MB) \cap \neg LA \cap \neg LB \Leftrightarrow x \in MA \cap \neg LA \cap \neg LB$ or $x \in MB \cap \neg LB \cap \neg LA \Leftrightarrow (x \in MA \cap \neg LA$ and $x \in \neg LB)$ or $(x \in MB \cap\neg LB$ and $x \in \neg LA)$$\Leftrightarrow (\mu_{A}(x) = 1/2$ and $\mu_{B}(x) \neq 1)$ or $(\mu_{B}(x) = 1/2$ and $\mu_{A}(x) \neq 1) \Leftrightarrow max(\mu_{A}(x), \mu_{B}(x)) =1/2$ or $max(\mu_{A}(x), \mu_{B}(x)) = 1/2 \Leftrightarrow max(\mu_{A}(x), \mu_{B}(x))=1/2$

$\lambda_{A \uplus B}(x) = 0 \Leftrightarrow x \in \neg M(A \uplus B) = \neg (MA \cup MB) = \neg MA \cap \neg MB \Leftrightarrow x \in \neg MA$ and $x \in \neg MB \Leftrightarrow \mu_{A}(x) = 0$ and $\mu_{B}(x) = 0 \Leftrightarrow max(\mu_{A}(x),\mu_{B}(x)) = 0$.

\smallskip

Villeurbanne, France, 1996

\

\small{
\noindent Added here to ease reading [12, p.160]:

\noindent $(\ast)$ $= (LU, MU) $ with $U = MX \cap Y \cap (X \cup M \neg Y) = (X \cap Y) \cup (MX \cap Y \cap M \neg Y)$

\noindent $(\ast \ast)$ $= (LV, MV)$ with $V = (X \cup Y) \cap (LX \cup Y \cup L \neg Y) = LX \cup Y \cup (X \cap L \neg Y)$}
\end{document}